\documentclass{amsproc}
\usepackage{amsfonts}
\usepackage{amssymb}
\usepackage{amscd}
\usepackage{verbatim}
\usepackage{graphicx}
\usepackage{subfigure}
\usepackage{fancyhdr}
\usepackage[english]{babel}
\usepackage[T1]{fontenc}
\usepackage[latin1]{inputenc}
\usepackage{bm}
\usepackage[colorlinks=true,urlcolor=blue]{hyperref}
\input xy
\xyoption{all}

\newcommand{\mor}[3]{$\xymatrix@1@C=15pt{#3: #1\ar[r]& #2}$}
\newcommand{\appl}[2]{$\xymatrix@1@C=15pt{#1 \ar@{|->}[r]& #2}$}
\newcommand{\pf}{\medskip \noindent {\bf Proof :\ \ }}
\newcommand{\wt}[1]{\widetilde{#1}}

\newcommand{\wh}[1]{\widehat{#1}}

\newcommand{\rad}[1]{{\rm rad}#1}
\newcommand{\BQI}{${\mathcal B}(Q,I)$}

\begin{document}
\bibliographystyle{plain}

\begin{abstract}
In this note, we investigate how different fundamental groups of
presentations of a  fixed algebra $A$  can be. For  finitely many
finitely  presented groups $G_i$, we construct  an algebra $A$
such that all  $G_i$ appear as fundamental groups of presentations
of $A$.
\end{abstract}

\title{Fundamental groups and presentations of algebras}
\author[J.~C.~Bustamante and D.~Castonguay]{ \begin{tabular}{ll}
Juan Carlos Bustamante & Diane Castonguay\\
           & \\
jc.bustamante@usherbrooke.ca& diane@inf.ufg.br\\
Universit\'e de Sherbrooke & Instituto de Inform\'atica \\
Bishop's University& Universidade Federal de Goi\'as \\
\end{tabular}\\
\footnote{Subject classification: 16G20, 16E40. Keywords and
phrases: Fundamental groups, bound quivers, presentations of
algebras. }}\address{Universit\'e de Sherbrooke, 2500 Boulevard de
l'Universit\'e, Sherbrooke, J1K 2R1, Qu\'ebec, Canada \newline
\indent Universidade Federal de Goi\'as, Instituto de Inform\'atica,
Bloco IMF I, Campus II, Samambaia - Caixa Postal 131,
CEP 74001-970, Goi\^ania, GO, Brasil}\maketitle

\thispagestyle{empty}

%
%

\section*{Introduction}

Let  $A$   be  a  basic,   connected,  finite  dimensional algebra
over  an algebraically  closed field  $k$. By  a  result of
Gabriel \cite{G79},  there exists a unique  finite connected
quiver $Q$ and a two-sided  ideal $I$ of the path algebra $kQ$,
such that $A \cong  kQ/I$. Such a pair $(Q, I)$ is called a bound
quiver. The morphism \mor{kQ}{A\simeq kQ/I}{\nu}, as well as
$(Q,I)$ are both called {\bf presentations} of  $A$. Following
\cite{MP83}, one can define the fundamental group  $\pi_1(Q,I)$.
Moreover, by a result  of Fischbacher and de  la Pe\~na
\cite{FdlP84}, every  finitely presented  group arises  in this
way. An important  feature of these groups is that  they depend
essentially on the ideal  $I$, thus it is  not an invariant  of
the algebra. Actually, there  are known examples of algebras $A
\cong kQ / I_1 \cong kQ/ I_2$ such that $\pi_1(Q, I_1) \ncong
\pi_1(Q, I_2)$ (see  example 1 in \ref{subsec:fund-groups}, and
section \ref{sec:change}).

With this  in mind,  we tackle  the following question:  Given an
algebra $A$, as above, how distinct can be fundamental groups of
presentations of $A$?

We  consider essentially  two  settings,  namely the  triangular
and the  not triangular case. In the triangular case, we consider
groups which are obtained from finite free  products of finitely
generated abelian  groups.  If we allow loops, we are able to
obtain results concerning finitely generated groups.

This  paper  is  organized  as  follows: In  Section
\ref{sec:Prel}, we  fix notations and terminology, recall
definitions concerning fundamental groups of bound quivers and
give some Examples. In Section \ref{sec:Prods}, we deal with
products  (and coproducts) of bound  quivers  which yield  to
products  (and coproducts) of their fundamental groups.  Section
\ref{sec:change} is devoted to investigate the  effects  of
changes of  presentations  on products  (and coproducts)  and
fundamental  groups. Finally,  in Section  \ref{sec:Teo}, we prove
the main result:

\subsection*{Theorem A}{\em Let $G_1, \ldots, G_n$ be finitely
presented groups. Then, there exists an algebra $A$ having
presentations $A\simeq kQ_A/I_i$, for $i \in \{1, \ldots, n\}$,
such that $\pi_1(Q_A,I_i)\simeq G_i$.}

The quivers considered in Theorem A have loops, so they lead to
algebras of infinite global dimension. Nevertheless, we obtain a
weaker result concerning triangular algebras. In this setting, we
consider the family of groups obtained from cyclic groups by
performing finite free and direct products, which we denote
$\mathbb{G}$.

\subsection*{Theorem B}{\em Let $G_1, \ldots, G_n \in \mathbb{G}$.
Then, there exists a triangular algebra $A$ having presentations
$A\simeq kQ_A/I_i$, for $i \in \{1, \ldots, n\}$, such that
$\pi_1(Q_A,I_i)\simeq G_i$.}

%
%

\section{Preliminaries}\label{sec:Prel}

\subsection{Bound quivers and algebras} A quiver $Q$ is a quadruple $(Q_0,
Q_1,$  $s,   t)$,  where  $Q_0$   and  $Q_1$  are   sets,  and
$s$,   $t$  maps \mor{Q_1}{Q_0}{s,t}.  The  elements of  $Q_0$ are
the  {\bf vertices}  of $Q$, whereas the elements of $Q_1$ are its
{\bf arrows}. Given an arrow $\alpha$ in $Q_1$, the vertex
$s(\alpha)$ is called its source and $t(\alpha)$ its target, and
we  write \mor{s(\alpha)}{t(\alpha)}{\alpha}.   A  {\bf  path} $w$
is  a sequence of arrows  $w=\alpha_1 \alpha_2\cdots \alpha_n$
such that  for $i \in \{1,\ldots,n-1\}$, one  has $t(\alpha_i)  =
s(\alpha_{i+1})$.  The  source and the terminus of a path $w$ in
$Q$ are defined in the obvious way. A quiver $Q$ is said to  be
{\bf finite} if both,  $Q_0$ and $Q_1$ are finite  sets. We say
that $Q$ is  {\bf connected}  if the  underlying graph  of $Q$ is
connected. Unless it is otherwise stated, we will consider finite
and connected quivers.

Given a commutative field  $k$ and a quiver $Q$, the path  algebra
$kQ$ is the $k-$vector space  whose basis is the set of paths of
$Q$, including  one stationary path  $e_x$ for  each vertex $x$ of
$Q$.   The multiplication  of two  basis elements of $kQ$ is their
composition whenever  it  is  possible, and $0$ otherwise. Let $F$
be the two-sided ideal of $kQ$ generated by the arrows of $Q$. A
two-sided ideal $I$ of  $kQ$ is called {\bf admissible} if there
exists an integer  $m\geq 2$  such that $F^m \subseteq I \subseteq
F^2$.   The pair $(Q,I)$ is then called a {\bf  bound quiver}.
Naturally, a {\bf pointed bound quiver} $(Q,I,x)$ is a bound
quiver together with a distinguished vertex $x\in Q_0$.

Conversely,  let $A$  be a  finite dimensional  algebra over  an
algebraically closed field  $k$. It is well-known  (see
\cite{G79,BG82} for example)  that, if in addition we  assume that
$A$ is a basic and  connected, then there  exists a unique finite
connected quiver $Q_A$ and a surjective  morphism of $k-$algebras
\mor{kQ_A}{A}{\nu}, which is  not unique, with $I={\rm Ker}\  \nu$
an admissible ideal. Those  surjective morphisms,  or equivalently
the  pairs $(Q_A,  I)$, are called  {\bf  presentations} of  the
algebra  $A$.   Remark that  a  morphism \mor{kQ_A}{A}{\nu}  is a
presentation of  $A$ whenever  $\{\nu (e_x)  \mid  x \in
(Q_A)_0\}$ is a complete  set of  primitive orthogonal idempotents
and,  for any fixed  $x,  y \in  (Q_A)_0$,  we  have  that $\{ \nu
(\alpha ) + \rad  ^2  A|\ $\mor{x}{y}{\alpha} $ \in Q_1\}$ is a
basis of $\nu (e_x)(\rad A/\rad ^2 A)\nu (e_y)$.

In this  note, for a given bound  quiver $(Q, I)$, we  will
consider morphisms \mor{kQ}{kQ/I \cong A}{\nu} defined by $\nu
(e_x)  = e_x + I$ for $x\in Q_0$, and, given an arrow $\alpha \in
Q_1$, from, say $x$ to $y$, $\nu (\alpha) = \alpha + \rho_\alpha +
I$ where $\rho _\alpha$ is a linear combination of paths from $x$
to $y$ different  from $\alpha$.   In  particular,  if  the paths
appearing  in $\rho_\alpha$ have length at least  2, then
\mor{kQ}{A}{\nu} is a presentation of $A$.

We refer the reader to \cite{ARS95, ASSk}, for instance, for
further reference on the use of bound quivers in the
representation theory of algebras.

\subsection{Fundamental groups of bound quivers}\label{subsec:fund-groups}
Given a bound quiver $(Q,I)$, its fundamental group is defined as
follows (see \cite{MP83}). For  $x,\ y$  in $Q_0$,  set $I(x,y) =
e_x (kQ)  e_y \cap  I$. A relation $\rho = \sum_{i=1}^m \lambda_i
w_i \in I(x,y)$ (where $\lambda_i \in k_*$,  and $w_i$  are
different paths from  $x$ to  $y$) is  said to  be {\bf minimal}
if $m\geq 2$, and, for every proper subset $J$ of $\{1, \ldots,
m\}$, we have $\sum_{i\in  J}  \lambda_i w_i  \notin  I(x,y)$. For
a given arrow \mor{x}{y}{\alpha}, let \mor{y}{x}{\alpha^{-1}} be
its formal inverse. A {\bf walk} $w$ in $Q$ from $x$  to $y$ is a
composition $w = \alpha_1^{\epsilon_1} \alpha_2^{\epsilon_2}
\cdots \alpha_n^{\epsilon_n   }$  such  that  the  $s(
\alpha_1^{\epsilon_1}) =  x,\ t(\alpha_n^{\epsilon_n}) =  y$, and,
for  $i \in \{2,\ldots n\}$, $s(\alpha_i^{\epsilon_i}) =
t(\alpha_{i-1}^{\epsilon_{i-1}})$.  Define the  {\bf homotopy
relation} $\sim$ on  the  set  of  walks  on  $(Q,I)$, as  the
smallest  equivalence  relation satisfying the following
conditions :

\begin{enumerate}

   \item For each arrow $\alpha$ from $x$ to $y$, one has $\alpha \alpha^{-1}
         \sim e_x$ and $\alpha^{-1} \alpha \sim e_y$.

   \item For each minimal relation $\sum_{i=1}^m \lambda_i w_i$, one has $w_i
         \sim w_j$ for all $i, j$ in $\{1, \ldots, m\}$.

   \item If $u,v, w$ and $w'$ are walks, and $u \sim v$ then $w u w' \sim w v
         w'$, whenever these compositions are defined.

\end{enumerate}

We denote by $\wt{w}$ the homotopy class  of a walk $w$.  Let
$v_0$ be a fixed point in $Q_0$, and consider the set $W(Q, v_0)$
of walks of source and target $v_0$. On this set, the product of
walks is everywhere defined. Because of the first and the  third
conditions in the definition of  the relation $\sim$, one can form
the  quotient group $W(Q, v_0)/\sim$.  This group  is called the
{\bf fundamental group} of the bound quiver $(Q, I)$ with base
point $v_0$, denoted by $\pi_1(Q,I,v_0)$. It follows easily from
the connectedness of $Q$ that this group does  not depend on  the
base  point $v_0$, and  we denote it  simply by $\pi_1(Q,I)$. This
group  has a clear geometrical interpretation  as the first
homotopy  group   of  a  C.W. complex  \BQI\  associated  to
$(Q,I)$,  see \cite{Bus04} (also \cite{K93}).

\subsection*{Remark} An important  remark, which is the main  motivation of
this work, is that the group defined above depends essentially  on
the minimal relations of the ideal $I$. It is well-known that, for
a $k-$algebra $A$, its presentation as a bound quiver algebra is
not unique. Thus,  the fundamental group is not an invariant of
the algebra, as the following well-known example shows (see also
section \ref{sec:change})

\subsection*{Example 1} Consider    the    quiver
$\xymatrix{3\ar@/^/[r]^\beta    \ar@/_/[r]_\gamma& 2\ar[r]^\alpha
&1}$  bound by the ideal $I_1 = <\beta  \alpha>$. Since the ideal
$I_1$ is generated  by  monomial relations,  the homotopy relation
is trivial,  thus $\pi_1(Q,I_1) \simeq \mathbb{Z}$. On the  other
hand, consider the morphism of algebras \mor{kQ}{A}{\nu_2} defined
by $\nu_2(\beta) = (\beta + \gamma)+ I_1$, and  $\nu_2(\alpha)  =
\alpha  +I_1$,  $\nu_2 (\gamma)  = \gamma  +I_1$.   A
straightforward computation shows that $\nu_2$ is a presentation,
and, moreover, that $I_2  = {\rm  Ker}\ \nu_2 = <(\beta - \gamma)
\alpha>$. This  yields to a trivial group $\pi_1(Q,I_2)$.
\smallskip

On the other  hand, these groups are invariant for  some classes
of algebras, as the following theorem states. Recall that a
$k-$algebra $A=kQ/I$ is said to be {\bf constricted} if for every
arrow \mor{x}{y}{\alpha} in $Q$ one has ${\rm dim}_k e_x  A e_y
=1$.

\subsection*{Theorem}{\bf (Bardzell - Marcos \cite{BM01})}{\em\ Let $A=kQ/I$ be
a constricted algebra. Given $(Q,I_1)$,  $(Q,I_2)$, two
presentations of $A$ one  then  has $\pi_1(Q,I_1) \simeq
\pi_1(Q,I_2)$.}\qed

\medskip It has  been proved in \cite{FdlP84} that  given a finitely presented
group $G$,  there exists an  incidence algebra $A=kQ/I$ such  that
$\pi_1(Q,I) \simeq G$.  A  natural question is then: {\em How
distinct can be fundamental groups of presentations of an
algebra?} Recalling that incidence algebras are always constricted
(and triangular),  and in light  of the Bardzell  - Marcos
theorem, this class of algebras is not interesting in view of our
problem.

The following example shows how given any finitely  presented
group $G$, one can  obtain  a  finite-dimensional  non constricted
nor  triangular  algebra $A\simeq kQ/I$, such that $\pi_1(Q,I)
\simeq G$.

\subsection*{Example 2} Let $G = <\alpha_i |\ w_j, 1\leq i \leq n, 1\leq j \leq
m>$ be a finitely presented group.  More precisely, $G$ is the
factor group of the free group  having basis
$\{\alpha_i\}_{i=1}^{n} $ by  the normal subgroup generated by
$\{w_i\}_{i=1}^{m}$. Without  loose of generality, we assume that
the  words $w_j$ are  reduced, non-empty,  pairwise different,
and, moreover, that each $w_j$ contains at least one letter among
the $\alpha_i$'s (otherwise we can replace $w_j$ by $w_j^{-1}$).
Consider the quiver $Q_G$: 
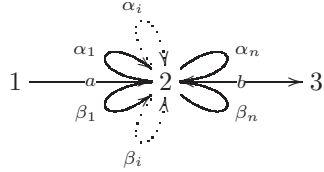
$$\xymatrix@C=20pt{1      \ar[rr]|a& &
2 \ar@(l,ul)[]^{\alpha_1} \ar@(ul,u)@{.>}[]^{\alpha_i}
\ar@(ur,r)[]^{\alpha_n}  \ar@(l,dl)[]_{\beta_1}
\ar@(dl,d)@{.>}[]_{\beta_i} \ar@(dr,r)[]_{\beta_n} \ar[rr]|b &&
3}$$ Identifying the arrow  $\beta_i$ with the letter
$\alpha_i^{-1}$,  we obtain a correspondence between the words  in
$\{\alpha_i, \alpha_i^{-1} \}_{i=1}^n$ and (some of)  the paths of
$Q_G$.   With this identification in  mind, define the ideal
$$I_G=<a\alpha_i\beta_ib - ab,\ aw_jb - ab, F^N|1\leq i\leq n,
1\leq j \leq m>$$ where $N={\rm  max}\{l(w_j)+3, 6 | 1\leq j \leq
m\}$.   We claim that the generators of  $I_G$ which  are not
monomial  relations are in  fact minimal relations.

Indeed, if this is not the case, then we have $ab \in I_G$. That is, there are
scalars  $\lambda_i,\ \mu_j$ with  $1\leq i  \leq n,  1\leq j  \leq m$  and an
element in $F^N$ such that  $$ab = \sum_{i=1}^n \lambda_i (a\alpha_1 \beta_i b
-ab)  +  \sum_{j=1}^m \mu_j(aw_jb  -ab  ) +  \gamma$$  \noindent  and this  is
equivalent  to  $$( \sum_{i=1}^n  \lambda_i  +  \sum_{j=1}^m  \mu_j -  1)ab  +
\sum_{i=1}^n  \lambda_i a\alpha_i \beta_i  b +  \sum_{j=1}^m \mu_j  a w_j  b +
\gamma = 0$$

\noindent But then, since $w_j \not= \alpha_i \beta_i$ for all $i,
j$, and $\gamma$ is a linear combination of paths  of length
$N\geq l(w_i) + 3 \geq 5$ we obtain $\lambda_i =  \mu_j = 0$ for
all $i,j$ and $\gamma  = 0$, which is a contradiction.

Therefore, $a\alpha_i \beta_i  b \sim ab$, and $aw_j b \sim  ab$.
From this we get $\alpha_i  \beta_i \sim  e_2$, and $w_j \sim
e_2$, for all $i,  j$. This shows that $\pi_1(Q_G,I_G)\simeq G$.
In section \ref{sec:change} we will consider changes of
presentations of $A=kQ_G/I_G$.

\medskip If one wishes to avoid quivers having oriented cycles, but
still have algebras that are not constricted, things are more
difficult.

\subsection{Coverings}\label{subsec:coverings}
Let $(\wh{Q}, \wh{I})$ be a (possibly  infinite) bound quiver, and
$G$ a group of automorphisms  acting freely on  $(\wh{Q},
\wh{I})$.  This ensures  that we can  form  the  quotient  $(Q,I)
=  (\wh{Q},  \wh{I})/G$.   The  natural  map \mor{(\wh{Q},
\wh{I})}{(Q,I)}{p} is  then called  a {\bf  covering}  of bound
quivers. In  this situation,  there is a  normal subgroup $H$  of
$\pi_1(Q,I)$ such that  $\pi_1(\wh{Q}, \wh{I}) \simeq  H$, and
$\pi_1(Q,I)/H \simeq  G$. In particular,  if $\pi_1(\wh{Q},
\wh{I})  =1$, then  $\pi_1(Q,I) \simeq  G$ (see \cite{P86},  for
instance).   Again,  in  this situation,  the  analogy  with
coverings of topological  spaces is clear.  It is shown in
\cite{Bus04} that $\wh{\mathcal{B}}  = \mathcal{B}(\wh{Q},
\wh{I})$ is a  regular covering space     of    $     \mathcal{B}
=     \mathcal{B}(Q,I)$     with    ${\rm
Cov}(\wh{\mathcal{B}}/\mathcal{B}) \simeq G$.

\subsection*{Example} For $n\geq 2$, let $\wh{Q}^n$ be the quiver whose
vertices are  $x_{i,j}$, with $0\leq i  < n$, $1\leq  j \leq n$
The  arrows of $\wh{Q}^n$        are
\mor{x_{i,j}}{x_{i,j-1}}{\alpha_{i,j}}        and
\mor{x_{i,j}}{x_{i+1,j-1}}{\beta_{i,j}} for  $0 \leq  i \leq n$,
$1\leq j\leq n$, where indices are to be read modulo $n$.
Moreover, let $\wh{I} = < w - w' | s(w)=s(w'),\  t(w)=  t(w')>$.
Consider  the  automorphism  \mor{(\wh{Q}^n, \wh{I})}{(\wh{Q}^n,
\wh{I})}{g}   defined  on   the  vertices  of   $Q^n$  by
$g(x_{i,j})=   x_{i+1,j}$,  and  on   the  arrows   by  $
g(\alpha_{i,j})  = \alpha_{i+1,j}$,  and  $g(\beta_{i,j})=
\beta_{i+1,j}$.   Then  $g$ has  order $n$. We can  form the
quotient $(Q^n,I) =  (\wh{Q}^n, \wh{I})/<g>$. The vertex set of
the quiver $Q^n$ is given  by: $Q_0^n = \{x_0, x_1, \ldots,
x_n\}$, and the arrows  are \mor{x_j}{x_{j-1}}{\alpha_j,\
\beta_j}, for $1\leq  j \leq n$. Moreover, $I\ =\ <\alpha_n \cdots
\alpha_1 - \beta_n \cdots \beta_1,\ \alpha_i \beta_{i-1} - \beta_i
\alpha_{i-1}|\ 1 < i \leq n >$. An immediate computation shows
that $\pi_1(\wh{Q}^n,  \wh{I}) = 1$ (see also  \cite{GR01,
Bus02}), so that $\pi_1(Q^n,I) \simeq \newline <g> \simeq
\mathbb{Z}_n$.

The quivers $\wh{Q}^2$, and $Q^2$ look as follows:

$$
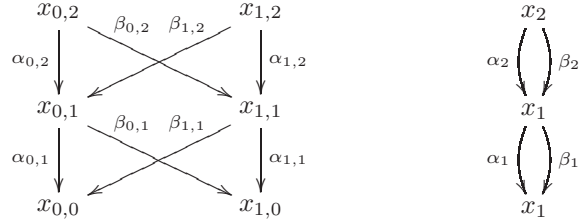
\xymatrix{   x_{0,2}  \ar[d]_{\alpha_{0,2}}   \ar[drr]^(.3){\beta_{0,2}}&  &
x_{1,2}   \ar[d]^{\alpha_{1,2}}   \ar[dll]_(.3){\beta_{1,2}}   &   &   &   x_2
\ar@/_/[d]_{\alpha_2} \ar@/^/[d]^{\beta_2}\\
x_{0,1}    \ar[d]_{\alpha_{0,1}}   \ar[drr]^(.3){\beta_{0,1}}&    &   x_{1,1}
\ar[d]^{\alpha_{1,1}} \ar[dll]_(.3){\beta_{1,1}}& &
       & x_1 \ar@/_/[d]_{\alpha_1} \ar@/^/[d]^{\beta_1}\\
        x_{0,0} & &x_{1,0} & &&x_1}$$

%
%
\section{Coproducts and Products}\label{sec:Prods}

Given  two groups, say  $G_1$ and  $G_2$, one  can consider  at least  two new
groups, namely  the free product $G_1 *  G_2 = G_1\amalg G_2$,  and the direct
product  $G_1  \times G_2$.   In  this  section, we  show  how  to carry  this
constructions to fundamental  groups of bound quivers.  Again,  the main ideas
come from algebraic topology. On one hand Van Kampen's theorem tells that some
fundamental groups  of topological spaces  are push-outs of groups,  thus free
products are involved. On the other hand, given two pointed topological spaces
$(X,x_0)$ and  $(Y,y_0)$ one can form  the product $(X  \times Y, (x_0,y_0))$,
and   then  $\pi_1(X   \times   Y,  (x_0,y_0))   \simeq  \pi_1(X,x_0)   \times
\pi_1(Y,y_0)$.

\subsection{Co-products} Given two pointed bound quivers $Q' = (Q', I', x')$
and $Q''= (Q'', I'',x'')$, assume, without loss of generality,
that $Q'_0 \cap Q_0'' = Q'_1 \cap Q_1'' = \emptyset$. We define
the quiver $Q = Q' \amalg Q''$ in the following way: $Q_0$ is
$Q'_0 \cup Q_0''$ in which we identify $x'$ and $x''$ to a single
new vertex $x$,  and $Q_1 = Q'_1 \cup Q''_1$. Then, $Q'$ and $Q''$
are identified  to two full convex sub-quivers of $Q$,  so walks
on $Q'$ or $Q''$  can be considered  as walks on  $Q$.  Thus, $I'$
and  $I''$ generate two-sided ideals of  $kQ$ which we denote
again by $I'$  and $I''$.  We define $I$ to  be the ideal  $I' +
I''$ of $kQ$. It  follows from this definition that the minimal
relations of $I'$ together with the minimal relations  of  $I''$
give the minimal relations needed to determine the homotopy
relation in $(Q,I)$. In addition, we can  consider an element
$\wt{w}\in \pi_1(Q',I',x')$  as an element $\wt{w} \in
\pi_1(Q,I,x)$. Conversely, any (reduced) walk $w \in W(Q,x)$ has a
decomposition $w=w'_1w''_1w'_2w''_2\cdots w'_nw''_n$ where $w'_i
\in W(Q',x')$, and   $w''_i\in  W(Q'',x'')$,  for
$i\in\{1,\ldots,n\}$, which is unique up to reduced walk.   In
addition, this decomposition is compatible with the homotopy
relations  involved. This leads us to the following proposition.

\subsection*{Proposition}\label{subsec:coproducts} With the notations above we
have: \begin{enumerate}

\item [$i)$] $(Q,I,x)$ is the coproduct, in the category of pointed bound
quivers, of $(Q', I',x)$ and $(Q'',I'', x)$

\item [$ii)$]  $\pi_1(Q,I,x) \simeq
\pi_1(Q',I',x') \amalg \pi_1(Q'',I'',x'')$.

\end{enumerate}

\pf The  first statement follows from  a direct computation,
while the second follows immediately from the above
discussion.\qed

\subsection*{Remark} It is worth to note that the canonical morphisms of pointed bound
quivers \mor{(Q',  I', x')}{(Q,  I, x)}{j'} and  \mor{(Q'', I'',
x'')}{(Q, I, x)}{j''} do  not induce morphisms of  $k-$algebras.
If one wants  to consider morphisms of  $k-$algebras, the  arrows
must be  {\it reversed}. One  then has canonical    projections
of    $k-$algebras    \mor{kQ'/I'}{kx}{p'}    and
\mor{kQ''/I''}{kx}{p''}. A natural question then is if the diagram

$$\xymatrix@C=50pt{kQ/I \ar[r] \ar[d]& kQ'/I'\ar[d]\\ kQ''/I''\ar[r]& kx}$$

\noindent is a pull-back of $k-$algebras. The answer is no. In
\cite{Lev04}      it     was     shown      that     the      pull-back     of
$\xymatrix@1@C=15pt{kQ'/I\ar[r]&kx & kQ''/I'' \ar[l]}$ is $kQ / J$ where $J= I
+ < Q_1'Q_1''> + <Q_1'' Q_1'>$.

\medskip
We now turn our interest into direct products.

\subsection{Products}\label{subsec:def-product} As before, consider two pointed bound
quivers $(Q',I',x')$, and $(Q'',I'', x'')$  whose source maps and targets maps
are  $s', s'',  t'$, and  $t''$. Following  \cite{Les94}, we  define  the {\bf
product quiver} $Q=Q' \otimes Q''$ as follows.  The vertex set $Q_0$ is simply
$Q'_0 \times  Q''_0$, whereas the arrow  set is $Q_1=  (Q'_1 \times Q''_0)\cup
(Q'_0   \times  Q''_1)$.   Given   an  arrow   $(\alpha,  y)   \in  Q'_1\times
Q''_0\subseteq Q_1$, define $s(\alpha, y ) = (s'(\alpha),y)$, and $t(\alpha,y)
=  (t'(\alpha),y)$.    In  an  analogous  way,  we   define  $s(x,\beta)$  and
$t(x,\beta)$  for  an  arrow  $(x,\beta)  \in Q'_0  \times  Q''_1$.   Now  let
$x=(x',x'')$ be the distinguished vertex in $Q$.

Given a  vertex $y  \in Q''_0$,  a path $w=\alpha_1  \cdots
\alpha_r$  in $Q'$ induces a  path $(\alpha_1, y)  \cdots \
(\alpha_r,y)$  in $Q$, which  we will denote by  $(w,y)$.
Similarly, for every  vertex $y$ of $Q'$, any path $w$ in $Q''$
induces a path $(y,w)$ in $Q$.

We define $I$ to be the  ideal whose generators are the following relations in
$kQ$:

\begin{enumerate}

\item [$(a)$] $(\rho, y)$, for every generator $\rho$ of $I'$, and
every vertex $y \in Q''$,

\item [$(b)$] $(x, \rho)$, for every generator $\rho$ of $I''$,
and every vertex $x \in Q'$,

\item [$(c)$] $(x_1,\beta)(\alpha,y_2) - (\alpha,y_1)(x_2, \beta)$
for every arrow \mor{x_1}{x_2}{\alpha} in $Q'$ and every arrow
\mor{y_1}{y_2}{\beta} in $Q''$.
\end{enumerate}

$$\xymatrix@C=50pt{(x_1,y_1)   \ar[r]^{(x_1,\beta)}   \ar[d]_{(\alpha,  y_1)}&
(x_1,y_2)\ar[d]^{(\alpha, y_2)}\\ (x_2,y_1)\ar[r]_{(x_2,\beta)}& (x_2,y_2)}$$

With these notations we have  an isomorphism of $k-$algebras $kQ'/I' \otimes_k
kQ''/I''  \simeq  kQ/I$ (see  \cite{Les94}).   However,  note  that since  the
natural projections  from $(Q,I,x)$ to $(Q', I',x')$  and $(Q'',I'',x'')$ {\em
are not}  morphisms of pointed  bound quivers, the  product quiver is  not the
product of $(Q',I',x')$ and $(Q'',  I'',x'')$ in the category of pointed bound
quivers.

Nevertheless,  given  an  arrow  $\Theta$  in  $Q$, we  can define
the  path $f'(\Theta)$ in $Q'$ in the following way

$$ f'(\Theta) = \left\{ \begin{array}{cl}
           \theta & \mbox{ if }  \Theta = (\theta,y) \in Q'_1\times Q''_0, \\
           e_x & \mbox{ if } \Theta = (x,\theta) \in Q'_0\times
Q''_1. \end{array} \right.$$

Extending  this map  in the  obvious way  to  walks in  $Q$, we
obtain a  map \mor{W(Q,x)}{W(Q',x')}{f'} which is, in fact, a
group homomorphism. In the same way, we obtain
\mor{W(Q,x)}{W(Q'',x'')}{f''}. This leads us to the following
lemma.

\subsection*{Lemma}{\em The maps $f'$ and $f''$ defined above induce groups
homomorphisms $\xymatrix@1@C=15pt{\phi':\pi_1(Q,I,x) \ar[r]& \pi_1(Q',I',x')}$
and $\xymatrix@1@C=15pt{\phi'':\pi_1(Q,I,x) \ar[r]& \pi_1(Q'',I'',x'')}$ given
by the rules $\phi'(\wt{w}) = \wt{f'(w)}$ and $\phi''(\wt{w}) = \wt{f''(w)}$}.

\pf The  only thing  one has  to show is  that $\phi'$  and
$\phi''$  are well defined maps.  We  do so only for $\phi'$. For
this  sake, define the relation $\approx$ in  $W(Q,x)$ by  $w_1
\approx  w_2$ if and  only if  $\wt{f'(w_1}) = \wt{f'(w_2)}$.
Keeping in  mind the generators of $I$, it  is easily seen that
this  is an  equivalence relation  on $W(Q,x)$  which verifies the
conditions $(1),\  (2),$ and  $(3)$  of the  definition  of the
homotopy  relation on  this set. Since  the latter is the smallest
such relation, we have  that $w_1 \sim w_2$   implies  $w_1
\approx  w_2$,   that   is  $w_1   \sim  w_2$   implies
$\phi'(\wt{w_1})  =  \phi'(\wt{w_2})$.   This   shows  that
$\phi'$  is  well defined.\qed

This leads us to the following proposition.

\subsection*{Proposition}\label{subsec:product}{\em With the notations above, we have an
isomorphism of  groups $$\pi_1(Q, I, x)  \simeq \pi_1(Q',I',x')
\times \pi_1(Q'',I'',x'').$$}

\pf In  light of  the preceding lemma,  we already  have a morphism  of groups
\mor{\pi_1(Q,I,x)}{\pi_1(Q',I',x')     \times    \pi_1(Q'',I'',x'')}{\Phi    =
(\phi',\phi'')}. In order to show that  this is an isomorphism, we exhibit its
inverse.

As noted before, given a walk $w$  in $Q'$, we can consider the
walk $(w,x'')$ in $Q$.  This  yields a map
\mor{\pi_1(Q',I',x')}{\pi_1(Q,I,x)}{\psi'} defined by
$\psi'(\wt{w)}  = \wt{(w,  x'')}$, which is, in fact a group
homomorphism. In the same way we obtain a group homomorphism
\mor{\pi_1(Q'',I'',x'')}{\pi_1(Q,I,x)}{\psi''}. This allows to
define  a group homomorphism \mor{ \pi_1(Q',I',x') \amalg
\pi_1(Q'',I'',x'')}{\pi_1(Q,I,x) }{\psi = \psi'  \amalg \psi''}.
Using relations of  type $c)$ in the definition of  the generators
of $I$,  one can easily  see that given  $\wt{w_1} \in
\pi_1(Q',I',x')$,  and $\wt{w_2}  \in \pi_1(Q'',I'',x'')$ one has
$\psi(\wt{w_1}\wt{w_2})       = \psi(\wt{w_2}\wt{w_1})$. Thus,
$\wt{w_1}\wt {w_2} \wt{w_1}^{-1}\wt{w_2}^{-1}$ belongs to ${\rm
Ker}\ \psi$, and, passing to the factor group, we can define a map
\mor{\pi_1(Q',I',x')\times\pi_1(Q'',I'',x'')}{\pi_1(Q,I)}{\Psi}
given by $\Psi(\wt{w_1},  \wt{w_2})  =   \psi(\wt{w_1}\wt{w_2})$.
Finally,  it  is  a straightforward   verification   that $\Phi$,
and  $\Psi$   are   mutually inverses.\qed

%
%

\section{Changes of presentations}\label{sec:change}
We  have  already  encountered  an  example  of  an  algebra  $A$
having  two presentations   $(Q,I_1)$  and  $(Q,I_2)$   such  that
$\pi_1(Q,I_1)  \simeq \mathbb{Z}$, and  $\pi_1(Q,I_2)$ is trivial.
Thus, we know how to  pass from an infinite cyclic group to a
trivial group. We begin this section by tackling the analogous
question for a finite (non-trivial) cyclic group $\mathbb{Z}_n$.

\subsection*{Example 1}\label{subsec:ex-changes-1} In section \ref{subsec:coverings} we considered the quiver $Q^n$:
$$\xymatrix{x_n\ar@/^/[r]^{\alpha_n}     \ar@/_/[r]_{\beta_n}
&    x_{n-1} \ar@/^/[r]^{\alpha_{n-1}}   \ar@/_/[r]_{\beta_{n-1}}
&   \cdots   &   \cdots \ar@/^/[r]^{\alpha_2} \ar@/_/[r]_{\beta_2}
&    x_1\ar@/^/[r]^{\alpha_1} \ar@/_/[r]_{\beta_1} & x_0}$$
\noindent bound by $I=<\alpha_n\cdots \alpha_1 - \beta_n\cdots
\beta_1,\ \alpha_i \beta_{i-1} - \beta_i \alpha_{i-1}|\ 1<i\leq n\
>$.  Recall that $\pi_1(Q^n,I) \simeq \mathbb{Z}_n$.

Consider the presentation  \mor{kQ^n}{kQ^n/I}{\nu} defined by
$\nu(\alpha_n) = (\alpha_n -  \beta_n) +  I$, and $\nu(\gamma)  =
\gamma  + I$ for  every arrow $\gamma \in Q^n_1,\  \gamma
\not=\alpha_n$. Let $J={\rm Ker}\  \nu$. We pretend that
$\pi_1(Q^n,J) \simeq 1$. In order to show this, let us begin by
computing the minimal relations of $J$.
   \begin{enumerate}
        \item Let $i\in \{2,\ldots,n\}$. Clearly, $\nu(\alpha_i \beta_{i-1} -
 \beta_i\alpha_{i-1})    \in    I$.    Thus,    $\alpha_i    \beta_{i-1}    -
 \beta_i\alpha_{i-1} \in J$, and this is a minimal relation.

       \item We pretend that $\alpha_n \beta_{n-1} + \beta_n\beta_{n-1} -
\beta_n \alpha_{n-1}$ is a minimal relation. To show this, it
suffices to show that its image by $\nu$ lies in $I$, since the
minimality is clear. We have: $$\begin{array}{l} \nu(\alpha_n
\beta_{n-1}  + \beta_n\beta_{n-1} -  \beta_n \alpha_{n-1})\\  =
\alpha_n  \beta_{n-1}  -  \beta_n \beta_{n-1}+ \beta_n\beta_{n-1}
-  \beta_n \alpha_{n-1} +I\\ = I
 \end{array}$$

       \item We pretend that $\alpha_n \cdots \alpha_1 + \beta_n \alpha_{n-1}
\cdots\alpha_1 - \beta_n \cdots \beta_1$ is a minimal relation.  Indeed,

$$\begin{array}{l}  \nu(\alpha_n   \cdots  \alpha_1  +   \beta_n  \alpha_{n-1}
\cdots\alpha_1 - \beta_n \cdots \beta_1)\\
= \alpha_n \alpha_{n-1} \cdots \alpha_1 - \beta_n \alpha_{n-1} \dots \alpha_1
+ \beta_n \alpha_{n-1}
\cdots\alpha_1 - \beta_n \cdots \beta_1+I\\ = I.
\end{array}$$

The relation 2,  above, gives us $\alpha_n \sim  \beta_n$, and
$\alpha_{n-1} \sim \beta_{n-1}$. Moreover, letting $i=n-1$  in
relations of type 2, we get $\alpha_{n-2} \sim  \beta_{n-2}$. The
same  argument, decreasing the  value of $i$ gives $\alpha_j \sim
\beta_j$ for $j\not=1$. Finally,  relation 3 gives $\alpha_1 \sim
\beta_1$. This shows that $\pi_1(Q^n,J)=1$.
\end{enumerate}

\medskip The following example illustrates how changes of presentations can be
done with  quivers having loops,  More precisely, it  shows how
one  can pass from any finitely presented group $G$ to the trivial
group.  \medskip

\subsection*{Example 2}\label{subsec:ex-changes-2} Let $G$ be a finitely presented group,
$(Q_G,I_G)$ the bound quiver of example 2, in section 1.2.
Moreover, set $A=kQ_G/I_G$. 
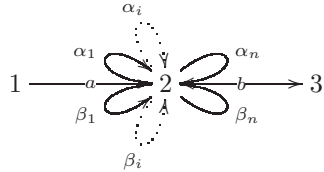
$$\xymatrix@C=20pt{1 \ar[rr]|a& & 2
\ar@(l,ul)[]^{\alpha_1} \ar@(ul,u)@{.>}[]^{\alpha_i}
\ar@(ur,r)[]^{\alpha_n} \ar@(l,dl)[]_{\beta_1}
\ar@(dl,d)@{.>}[]_{\beta_i} \ar@(dr,r)[]_{\beta_n} \ar[rr]|b &&
3}$$ \noindent Recall that $I_G =<a\alpha_i\beta_ib - ab,\ aw_jb -
ab, F^N|1\leq i\leq n, 1\leq j \leq m \}$ where  $N={\rm
max}\{l(w_j)+3,6| 1\leq  j \leq m\}$, and this leads to
$\pi_1(Q_G, I_G)\simeq G$.

Consider now the morphism \mor{kQ_G}{kQ_G/I_G}{\nu} defined on the arrows of
$Q_G$ by $$\nu(\gamma) = \left\{ \begin{array}{ll}
\alpha_i + \beta_i +\beta_i^2 + I_G& \hbox{ if } \gamma = \alpha_i, \\
\gamma + I_G & \hbox{ otherwise. } \end{array} \right.$$

\noindent Since $\{ \nu(\alpha_i) + {\rm rad}^2 A, \nu(\beta_i) +
{\rm rad }^2 A | 1\leq i \leq n \}  = \{\alpha_i + \beta_i + {\rm
rad} ^2 A, \beta_i + {\rm rad}^2 A \}$, this set is a basis of
$e_2 ({\rm rad} A / {\rm rad}^2 A) e_2$, so that $\nu$ is a
presentation of $A$. In particular ${\rm Ker}\ \nu$ is an
admissible ideal. We claim that $\pi_1( Q_G, {\rm Ker}\ \nu )
\simeq 1$.

We show that $\rho = a\alpha_i \beta_i b -a\beta_i^2b - a\beta_i^3
b-ab$ is a minimal relation. Indeed
\begin{eqnarray*}
\nu(\rho) & = & \nu( a\alpha_i \beta_i b -a\beta_i^2b - a\beta_i^3)\\
& = & a(\alpha_i + \beta_i + \beta_i^2)\beta_i b - a\beta_i^2b -a\beta_i^3 b
-ab + I_G\\
&=& a\alpha_i \beta_i b - ab +I_G\\
&=& I_G
\end{eqnarray*}
\noindent Thus, it only remains to prove the minimality of $\rho$,
but this follows from the following facts :

\begin{enumerate}

\item $\nu(ab) = ab + I_G \not=I_G$,

\item A linear combination of $a\beta_i^2 b$ and $a\beta_i^3 b$
cannot belong to $I_G$. Indeed, it follows from the hypothesis
made on the words $w_i$ that each one of them contains at least
one $\alpha_i$, thus is different from $\beta_i^2$ and
$\beta_i^3$. In addition we have $$l(a\beta_i^2 b) < l(a\beta_i^3
b) = 5 < N.$$\end{enumerate}

\smallskip
\noindent So we are done, and $\rho$ is minimal. This yields$$
a\alpha_i \beta_i b \sim a\beta_i^2 b \sim a\beta_i^3 \sim ab$$ so
that $\alpha_i \beta_i \sim e_2$, and $\beta_i^2 \sim \beta_i^3$,
and this shows our claim.

\medskip
The encountered examples all show that there are algebras having
different presentations which have an arbitrary group, as well as
the trivial group as fundamental groups. The following example
shows how one can pass directly from any finitely generated
abelian group to a free abelian group.
\medskip

\subsection*{Example 3}\label{subsec:ex-changes-3}  Fix a positive integer $t$,
and, for each $i\in \{1,\ldots,t\}$, let $n_i\in \mathbb{N}_*$,
and let $m>{\rm max}_i \{n_i\} + 1$.  Consider the quiver
$$
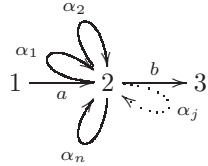
\xymatrix{1 \ar[r]_{a} & 2 \ar@(l,ul)[]^{\alpha_1}
\ar@(ul,u)[]^{\alpha_2} \ar@(r,dr)@{.>}[]^{\alpha_j}
\ar@(d,dl)[]^{\alpha_n}\ar[r]^b& 3}$$

\noindent bound by $I=<a \alpha_i b + a \alpha_i^{n_i + 1} b,\
\alpha_i \alpha_j - \alpha_j \alpha_i, F^m|\ 1\leq i \leq t
>$. A direct computation shows that $\pi_1(Q,I) \simeq
\bigoplus_{i=1}^t \mathbb{Z}_{n_i}$. On the other hand, consider
the presentation \mor{kQ}{kQ/I}{\nu} given by $\nu(\alpha_i) =
(\alpha_i - \alpha_i^{n_i+1}) +I$. We leave the reader verify that
$I'={\rm Ker}\ \nu = <a\alpha_i b, \alpha_i \alpha_j - \alpha_j
\alpha_i, F^m|\ 1\leq i,j\leq t>$, and this leads to $\pi_1(Q,I')
= \bigoplus_{i=1}^t \mathbb{Z}$.

Note that choosing $n_i=1$ for all $i\in\{1,\ldots,t\}$ we obtain
that $\pi_1(Q,I)$ is the trivial group.

\medskip
In section \ref{sec:Prods}, we saw that fundamental groups of
bound quivers behave well under products and coproducts. The
following lemma shows that the same is true under changes of
presentations.

\subsection*{Lemma}{\em  For $i\in \{1,2\}$, let $A_i \simeq kQ_{A_i}/I_i \simeq kQ_{A_i}/I_i'$
be algebras with two different presentations. Denote by
$\pi_1(Q_{A_i}, I_i) \simeq G_i$, and $\pi_1(Q_{A_i}, I'_i) \simeq
G'_i$. Then we have the following:

\begin{enumerate}

   \item [$i)$] There exists an algebra $A$ having two presentations
   $A\simeq kQ/I \simeq kQ/I'$ such that $\pi_1(Q,I) \simeq G_1
   \amalg G_2$ and $\pi_1(Q,I') \simeq G'_1
   \amalg G'_2$.\\

   \item [$ii)$] There exists an algebra $A$ having two presentations
   $A\simeq kQ/I \simeq kQ/I'$ such that $\pi_1(Q,I) \simeq G_1
   \times G_2$ and $\pi_1(Q,I') \simeq G'_1
   \times G'_2$.
\end{enumerate}}

\pf For $i\in \{1,2\}$, let \mor{kQ_{A_i}}{A_i}{\nu_i}, and
\mor{kQ_{A_i}}{A_i}{\nu'_i} be the presentations of the algebra
$A_i$ such that $I_i={\rm Ker}\ \nu_i$, and $I'_i={\rm Ker}\
\nu'_i$. In order to prove the first statement, consider the
pointed bound quiver $(Q,I,x) = (Q_{A_1}, I_1, x_1) \amalg
(Q_{A_2},I_2,x_2)$ and let $A=kQ/I$. It follows from theorem
\ref{subsec:coproducts} that $\pi_1(Q,I,x)\simeq G_1 \amalg G_2$.
On the other hand, consider the presentation \mor{kQ}{kQ/I\simeq A
}{\nu'} given by $\nu'(\alpha) = \nu'_i(\alpha)$, where $\alpha$
is an arrow of $Q_{A_i}$. Then, we have $I'={\rm Ker}\ \nu' = I'_1
+ I'_2$, and, again, from theorem \ref{subsec:coproducts}, we
obtain $\pi_1(Q,I',x)\simeq G'_1 \amalg G'_2$.

In order to prove the second statement, consider the quiver $Q =
Q_{A_1} \otimes Q_{A_2}$ bound by the ideal $I$ as described in
section \ref{subsec:def-product}. We then have, from theorem
\ref{subsec:product}, that $\pi_1(Q, I, x)\simeq G_1\times G_2$.
Let $A=kQ/I$, and consider the following presentation
\mor{kQ}{kQ/I \simeq A}{\nu'} given by

$$ \nu'(\Theta) = \left\{ \begin{array}{cl}
          \nu'_1( \theta )    & \mbox{ if } \Theta = (\theta,y) \in (Q_{A_1})_1 \times (Q_{A_2})_0, \\
          \nu'_2(\theta  )    & \mbox{ if } \Theta = (x,\theta) \in (Q_{A_1})_0 \times (Q_{A_2})_1. \end{array} \right.$$

Again, it follows from the definition of the ideal $I$, and using
the fact that $\nu'_1$ and $\nu'_2$ are changes of presentations,
that $I = {\rm Ker}\ \nu'$. Thus, theorem \ref{subsec:product}
gives  $\pi_1(Q, I', x)\simeq G'_1\times G'_2$.\qed

We are now able to prove our main results.

%
%

\section{Main Results}\label{sec:Teo}

\subsection*{Theorem A}{\em Let $G_1, \ldots, G_n $ be finitely presented groups.
Then, there exists a finite dimensional algebra $A$ having
presentations $A\simeq kQ_A/I_i$, for $i \in \{1, \ldots, n\}$,
such that $\pi_1(Q_A,I_i)\simeq G_i$.}

\pf  Using example 3.2, we can build algebras  $A_i$ having
presentations $A_i\simeq kQ_i/J_i \simeq kQ_i/J'_i$ with
$\pi_1(Q_i,J_i) \simeq G_i$, and $\pi_1(Q_i,J'_i)\simeq 1$. For
$i\in\{1,\ldots,n\}$ consider the bound quiver

$$(Q_A,I_i)  =\left( \coprod_{l=1}^{i-1}  (Q_l,J'_l)\right) \amalg
(Q_i, J_i) \amalg \left( \coprod_{l=i+1}^n(Q_l,J'_l)\right).$$
\noindent It follows from theorem \ref{subsec:coproducts},
statement $ii)$ that $\pi_1(Q_A,I_i) \simeq G_i$. Moreover, using
the argument of the proof of the above lemma, one gets $kQ_A/I_i
\simeq kQ_A/I_j$, for all $i,j\in\{1,\ldots,n\}$.\qed

\medskip
The examples of  quivers without oriented cycles that  we have
encountered led us   to  consider   cyclic  groups.    Moreover,
from   results   in  section \ref{sec:Prods}  we  know  how  to
deal  with  products  and  coproducts  of fundamental groups.  Let
us denote by $\mathbb{G}$ the smaller family of groups satisfying
the following conditions :

\begin{enumerate}
   \item If $G$ is a cyclic group then $G \in \mathbb{G}$,
   \item If $G_1,\ G_2$ belong to $ \mathbb{G}$, then the same
   holds for $G_1\coprod G_2$ and $G_1 \times G_2$.
\end{enumerate}

This leads us to the following theorem.

\subsection*{Theorem B}{\em Let $G_1, \ldots, G_n \in \mathbb{G}$.
Then, there exists a triangular algebra $A$ having presentations
$A\simeq kQ_A/I_i$, for $i \in \{1, \ldots, n\}$, such that
$\pi_1(Q_A,I_i)\simeq G_i$.}

\pf Using  propositions  \ref{subsec:coproducts},
\ref{subsec:product}, example 1 in section \ref{sec:change},
example 1 in \ref{subsec:fund-groups}, and the above lemma, we can
build triangular algebras $A_i$ having presentations $A_i\simeq
kQ_i/J_i \simeq kQ_i/J'_i$ with $\pi_1(Q_i,J_i) \simeq G_i$, and
$\pi_1(Q_i,J'_i)\simeq 1$. The remaining part of the proof is just
as in the preceding result.\qed

\subsection*{Corollary}{Let $M_1, \ldots, M_n$ be finitely
generated abelian groups. Then there exists a triangular algebra
$A$ having presentations $A\simeq kQ_A/I_i$, for $i \in \{1,
\ldots, n\}$, such that $\pi_1(Q_A,I_i)\simeq M_i$.\qed

\medskip
As a final remark, let us note that the homotopy relation in a
bound quiver does make sense even if we do not ask the ideal $I$
to be admissible, nor the quiver $Q$ to be finite. These
requirements lead to finite dimensional algebras (with $1$).

If one wishes to consider the family of all finitely generated
groups, similar constructions of what have been made in this work
can be done. Indeed given any finitely generated group $G$, as in
Example 2 of \ref{subsec:fund-groups} one can construct a quiver
$Q_G$ and a two sided ideal $I_G \trianglelefteq kQ_G$ such that
$\pi_1(Q_G, I_G)\simeq G$. In this case, the quiver would still be
finite, but the ideal $I$ will not be admissible, so this would
lead to infinite dimensional algebras which can be seen as locally
bounded $k-$categories.

As further generalisation, one may wish to consider arbitrary
groups. Again, this can be performed, as in Example 2 of
\ref{subsec:fund-groups}. This time the quiver will still have
three vertices, but will have an infinite number of arrows, and,
again, the ideal would not be admissible.

\medskip
Also, in light of Bardzell-Marcos theorem, one may ask how big is
the family of groups that rise as fundamental groups of
presentation of a fixed algebra $A$.

\subsection*{Acknowledgements} Both authors thank Pierre-Yves Leduc
for the interesting question which  brought us to this note. This
work was done while the second author was professor at the
University of S\~{a}o Paulo, and the first author had
post-doctoral fellowship at the same University. The first author
gratefully acknowledges the S\~{a}o Paulo group for hospitality
during his stay there, as well as financial support from
F.A.P.E.S.P., Brazil

\bibliography{biblio}
\end{document}